\documentclass{article}

\usepackage{amssymb}
\usepackage{latexsym}

\usepackage{graphicx}

\title{An extension of Motzkin-Straus Theorem to non-uniform hypergraphs and its applications}

\author
{Yuejian Peng \thanks{
College of Mathematics, Hunan University, Changsha 410082, P.R. China  Email: ypeng1@163.com.  Supported in part by National Natural Science Foundation of China (No. 11271116). }
 \and Hao Peng \thanks{College of Mathematics, Hunan University, Changsha 410082, P.R. China.Email: hpeng@hnu.edu.cn }
 \and Qingsong Tang \thanks{College of Sciences, Northeastern University, Shenyang, 110819, China and Mathematics School, Institute of Jilin University, Changchun, 130012, China.Email: t\_qsong@sina.com.cn}
 \and Cheng Zhao \thanks{Department of Mathematics and Computer Science, Indiana State University, Terre Haute, IN, 47809 and School of Mathematics, Jilin University, Changchun 130012, P.R. China. Email: cheng.zhao@indstate.edu} }

\date{}

\newtheorem{defi}{Definition}[section]
\newtheorem{theo}{Theorem}[section]

\newtheorem{remark}[theo]{Remark}

\newtheorem{lemma}[theo]{Lemma}
\newtheorem{claim}[theo]{Claim}

\newcommand{\qed}{\hspace*{\fill} \rule{7pt}{7pt}}

\begin{document}
\maketitle
\begin{abstract}
 In 1965,  Motzkin and Straus established a remarkable connection between the order of a  maximum clique and the Lagrangian of a graph and  provided a new proof of Tur\'an's theorem using the connection.  The connection of Lagrangians and Tur\'{a}n densities can be also used to prove the fundamental theorem of Erd\H os-Stone-Simonovits on Tur\'{a}n densities of graphs.  Very recently, the study of Tur\'{a}n densities of non-uniform hypergraphs have been motivated by extremal poset problems. In this paper, we attempt to explore the applications of Lagrangian method in determining Tur\'{a}n densities of non-uniform hypergraphs. We first give a  definition of the Lagrangian of a non-uniform hypergraph, then  give an extension of Motzkin-Straus theorem to non-uniform hypergraphs whose edges contain 1 or 2 vertices. Applying it, we give an extension of  Erd\H os-Stone-Simonovits theorem to non-uniform hypergraphs whose edges contain 1 or 2 vertices.
\end{abstract}


\section{Introduction and main results}

 Tur\'{a}n problems on uniform hypergraphs have been actively studied.   In 1965,  Motzkin and Straus provided a new proof of Tur\'an's theorem based on a remarkable connection between the order of a  maximum clique and the Lagrangian of a graph in \cite{MS}.  In fact, the connection of Lagrangians and Tur\'{a}n densities can be used to give another proof of the fundamental theorem of Erd\H os-Stone-Simonovits on Tur\'{a}n densities of graphs in \cite{keevash}. This type of connection aroused  interests in the study of Lagrangians of uniform hypergraphs. Very recently, the study of Tur\'{a}n densities of non-uniform hypergraphs have been motivated by extremal poset problems (see \cite{GK08} and \cite{GL09}). In this paper, we attempt to explore the applications of Lagrangian method in determining Tur\'{a}n densities of non-uniform hypergraphs. We first give a  definition of the Lagrangian of a non-uniform hypergraph, then  give an extension of Motzkin-Straus theorem to non-uniform hypergraphs whose edges contain 1 or 2 vertices. Applying it, we give an extension of  Erd\H os-Stone-Simonovits theorem to non-uniform hypergraphs whose edges contain 1 or 2 vertices.

A hypergraph $H=(V,E)$ consists of a vertex set $V$ and an edge set $E$, where every edge in $E$ is a subset of $V$. The set $R(H)=\{|F|:F \in E\}$ is called the set of  $edge\ types$ of $H$. We also say that $H$ is a $R(H)$-graph. For example, if $R(H)=\{1, 2\}$, then we say that $H$ is a $\{1, 2\}$-graph. If all edges have the same cardinality $k$, then $H$ is a $k$-uniform hypergraph. A $2$-uniform graph is called a graph. A hypergraph is non-uniform if it has at least two edge types.  For any $k \in R(H)$, the $level \ hypergraph \ H^k$ is the hypergraph consisting of all edges with $k$ vertices of $H$. We write $H^R_n$ for a hypergraph $H$ on $n$ vertices with $R(H)= R$. An edge $\{i_1, i_2, \cdots, i_k\}$ in a hypergraph  is simply written as  $i_1i_2\cdots i_k$ throughout the paper.

 The complete hypergraph $K^R_n$ is a hypergraph on $n$ vertices with edge set $ \bigcup_{i\in R} {[n] \choose i}$.   For example, $K_n^{\{k\}}$ is the complete k-uniform hypergraph on $n$ vertices. Let $[k]$ denote the set $\{1, 2, \cdots, k\}$, then  $K_n^{[k]}$ is the non-uniform hypergraph with all possible edges of cardinality at most $k$. The complete graph on $n$ vertices $K_n^{\{2\}}$ is also called a clique. We also let $[k]^{(r)}$ represent the complete $r$-uniform hypergraph on vertex set $[k]$.

  Let us briefly review  the Tur\'{a}n problem  on uniform hypergraphs. For a given $r$-uniform graph $F$ and positive integer $n$,
let ${\rm ex}(n, F)$ be  the maximum number of edges an $r$-uniform
graph on $n$ vertices can have without containing  $F$ as a subgraph. By
a standard averaging argument of Katona, Nemetz, and Simonovits in
\cite{KNS}, ${{\rm ex}(n,  F) \over {n \choose r}}$
decreases as $n$ increases, therefore 
$\lim_{n\rightarrow\infty}{{\rm ex}(n, F) \over {n
\choose r}}$ exists.  This limit is  called the {\it Tur\'an
density} of $F$ and denoted by $\pi(F)$.
Tur\'an's theorem
\cite{Turan} says that $\pi(K^{\{2\}}_{l})=1-{1 \over
l-1}$. The fundamental result in extremal graph theory due to  Erd\H{o}s-Stone-Simonovits generalizes Tur\'an's theorem and it says that for a graph $F$ with chromatic number $\chi(F)$ where $\chi(F) \ge 3$,
then $\pi(F)=1-{1 \over \chi(F)-1}$. However, we know quite few about Tur\'an density of $r$-uniform hypergraphs for $r\ge 3$ though some progress has been made.

A useful tool in extremal problems of uniform hypergraphs (graphs) is the Lagrangian of a uniform hypergraph (graph).

\begin{defi} \label{1.1}
Let $G$ be  an $r$-uniform graph  with vertex set $[n]$ and
edge set $E(G)$.  Let $S=\{\vec{x}=(x_1,x_2,\ldots ,x_n)\in R^n: \sum_{i=1}^{n} x_i =1, x_i
\ge 0 {\rm \ for \ } i=1,2,\ldots , n \}$. For $\vec{x}=(x_1,x_2,\ldots ,x_n)\in S$,
define
$$\lambda (G,\vec{x})=\sum_{i_1i_2 \cdots i_r \in E(G)}x_{i_1}x_{i_2}\ldots x_{i_r}.$$

 The Lagrangian of
$G$, denoted by $\lambda (G)$, is defined as $$\lambda (G) = \max \{\lambda (G, \vec{x}): \vec{x} \in S \}.$$

\end{defi}

Motzkin and Straus in \cite{MS} shows that  the Lagrangian of a graph is determined by the order of its maximum clique.

\begin{theo} (Motzkin and Straus \cite{MS}) \label{MStheo}
If $G$ is a graph in which a largest clique has order $l$, then
$\lambda(G)=\lambda(K^{\{2\}}_l)=\lambda([l]^{(2)})={1 \over 2}(1 - {1 \over l})$.
\end{theo}

This connection provided another proof of Tur\'an's theorem. More generally, the connection of Lagrangians and Tur\'{a}n densities can be used to give another proof of  Erd\H os-Stone-Simonovits result(see Keevash's survey paper \cite{keevash}). In 1980's,  Sidorenko \cite{sidorenko87} and Frankl and F\"uredi \cite{FF88} developed the method of  applying  Lagrangians in determining hypergraph Tur\'an densities. More applications of Lagrangians can be found in \cite{FF}, \cite{FR84}, \cite{FR89}, \cite{mubayi06}, \cite{sidorenko89}. Very recently, the study of Tur\'{a}n densities of non-uniform hypergraphs have been motivated by the study of extremal poset problems \cite{GK08}, \cite{GL09}. A generalization of the concept of Tur\'{a}n density to a non-uniform hypergraph was given in \cite{JL}.

For a non-uniform hypergraph $G$ on $n$ vertices, the Lubell function of $G$ is defined to be
$$h_n(G)=\sum_{k\in R(G)}{|E(G^k)| \over {n \choose k} }.$$

Given a family of hypergraph $\mathcal{F}$ with common set of edge-types $R$, the Tur\'{a}n density of $\mathcal{F}$ is defined to be
 \begin{eqnarray*}
 \pi(\mathcal{F})=\lim_{n\rightarrow\infty}\max\{h_n(G): |v(G)|=n, G \subseteq K^R_n, {\rm \ and\ } G {\rm \ contains \ no \ subgraph \ in \ } \mathcal{F} \}.
 \end{eqnarray*}
The proof of the existence of this limit can be found in \cite{JL}.
\begin{defi}
For any hypergraph $H_n$ and positive integers $s_1,s_2,\ldots,s_n$, the blowup of $H$ is a new hypergraph $(V,E)$, denoted by $H(s_1,s_2,\ldots,s_n)$, satisfying

$1. \,\  V=\bigcup_{i=1}^n{V_i},  \,\ where \,\ |V_i|=s_i;$

$2. \,\  E=\bigcup_{F\in E(H)}\Pi_{i\in F}V_i.$\\
\end{defi}

\begin{remark}\label{re2} For a non-uniform hypergraph $G$ on $n$ vertices, the blowup of $G$ has the following property:
$$ \lim_{t\rightarrow\infty}h_{nt}(G(t,t,\ldots,t))= h_n(G).$$
\end{remark}

The Lagrangian of a $k$-uniform graph $G$ is the supremum of the  densities of blowups of $G$ multiplying the constant ${1 \over k!}$ (see \cite{keevash}). We define the Lagrangian of a non-uniform hypergraph as follows so that the Lagrangian of a non-uniform hypergraph $H$ is the supremum of the densities of blowups of $H$.

\begin{defi}
For  a hypergraph $H^{R}_n$ with $R(H)=R$ and a vector $\vec{x}=(x_1,\ldots,x_n) \in R^n$,
define
$$\lambda' (H^{R}_n,\vec{x})=\sum_{j\in R}(j!\sum_{i_1i_2 \cdots i_j \in H^j}x_{i_1}x_{i_2}\ldots x_{i_j}).$$
\end{defi}

\begin{defi} \label{1.4}
Let $S=\{\vec{x}=(x_1,x_2,\ldots ,x_n): \sum_{i=1}^{n} x_i =1, x_i
\ge 0 {\rm \ for \ } i=1,2,\ldots , n \}$. The Lagrangian of
$H^{R}_n$, denoted by $\lambda' (H^{R}_n)$, is defined as
 $$\lambda' (H^{R}_n) = \max \{\lambda' (H^{R}, \vec{x}): \vec{x} \in S \}.$$
The value $x_i$ is called the {\em weight} of the vertex $i$. We call $\vec{x}=(x_1, x_2, \ldots, x_n) \in R^n$ a legal weighting for $H$ if $\vec{x} \in S$. A vector $\vec{y}\in S$ is called an {\em optimal weighting} for $H$ if $\lambda' (H, \vec{y})=\lambda'(H)$.
\end{defi}

\begin{remark}\label{relation} The connection between Definitions \ref{1.1} and {1.4} is that, 
if $G$ is a $k$-uniform graph, then $$\lambda'(G)=k!\lambda(G).$$
\end{remark}

In this paper, we will prove  the following generalization of Motzkin-Straus result to $\{1, 2\}$-graphs.

\begin{theo} \label{th1}
If $H$ is a $\{1,2\}$-graph and the order of its maximum complete $\{1,2\}$-subgraph is $t$ ( where $t\ge 2$), then
$\lambda'(H)=\lambda'(K^{\{1,2\}}_t)=2 - {1 \over t}$.
\end{theo}

As an application of Theorem \ref{th1}, we will also prove an extension of Erd\H{o}s-Stone-Simonovits result to $\{1, 2\}$-graphs as given in the following theorem. This result was proved by Johnston and Lu in \cite{JL} using a different approach. Our motivation is to explore the applications of Lagrangian method in the Tur\'an problem.

\begin{theo} \label{th2}
If  $H$ is a $\{1,2\}$-graph and $H^2$ is not bipartite,  then $\pi(H)=2-{1\over {\chi(H^2)-1}}$.
\end{theo}

\section{Proofs of the main results}

We will impose an additional condition on any optimal weighting ${\vec x}=(x_1, x_2, \ldots, x_n)$ for a hypergraph $H$:

(*) $ \vert\{i : x_i > 0 \}\vert$ is minimal, i.e., if $\vec y $ is a legal weighting for $H$ satisfying $|{i : y_i > 0}| < |{i : x_i > 0}|$, then $ \lambda' (G, {\vec y}) < \lambda'(G)$.

\subsection{Proof of Theorem \ref{th1}}

We need the following two lemmas.

\begin{lemma}\label{Lemma1}
If $x_1 \ge x_2 \ge \ldots \ge x_k >x_{k+1}=x_{k+2}=\ldots x_n=0$ and ${\vec x}=(x_1, x_2, \ldots, x_n)$ is an optimal weighting of a hypergraph $H$, then ${\partial\lambda' (H, {\vec x})\over \partial x_1}={\partial \lambda' (H, {\vec x})\over \partial x_2}= \ldots ={\partial \lambda' (H, {\vec x})\over \partial x_k}$ .
\end{lemma}

\noindent{\em  Proof.} Suppose, for a contradiction, that there exist $i$ and $j$ $(1\le i<j\le k)$ such that ${\partial\lambda' (H, {\vec x})\over \partial x_i}> {\partial \lambda' (H, {\vec x})\over \partial x_j}$. We define a new legal weighting $\vec y$ for $H$ as follows. Let $y_l=x_l$ for $l \neq i,j$, $y_i=x_i+ \delta$ and $y_j=x_j- \delta \geq 0,$ then
\begin{eqnarray*}
&&\lambda'(H, \vec y)- \lambda'(H, \vec x)\\
&=&\delta({\partial \lambda' (H, {\vec x})\over \partial x_i}- x_j{\partial^2 \lambda' (H, {\vec x})\over {\partial x_i \partial x_j}})-
\delta({\partial \lambda' (H, {\vec x})\over \partial x_j}- x_i{\partial^2 \lambda' (H, {\vec x})\over {\partial x_i \partial x_j}})+
(\delta x_j-\delta x_i-\delta^2){\partial^2 \lambda' (H, {\vec x})\over {\partial x_i \partial x_j}}\\
&=&\delta({\partial \lambda' (H, {\vec x})\over \partial x_i}- {\partial \lambda' (H, {\vec x})\over \partial x_j})-\delta^2{\partial^2 \lambda' (H, {\vec x})\over {\partial x_i \partial x_j}}\\
&>&0
\end{eqnarray*} 
for some small enough $\delta$, contradicting to that $\vec x$ is an optimal vector. Hence Lemma \ref{Lemma1} holds.\qed

\begin{lemma}\label{Lemma2}
If $x_1 \ge x_2 \ge \ldots \ge x_k >x_{k+1}=x_{k+2}=\ldots x_n=0$ and ${\vec x}=(x_1, x_2, \ldots, x_n)$ is an optimal weighting of a hypergraph $H$ satisfying (*), then $\forall 1\le i<j\le k,$ there exists an edge $e \in E(H)$  such that $\{i,j\} \subseteq e$.
\end{lemma}

\noindent{\em  Proof.} Suppose, for a contradiction, that there exist $i$ and $j$ $(1\le i<j\le k)$ such that $\{i,j\} \nsubseteq e$ for any $e \in E(H)$. We define a new weighting $\vec y$ for $H$ as follows. Let $y_l=x_l$ for $l \neq i,j$, $y_i=x_i+ x_j$ and $y_j=x_j- x_j=0,$ then $\vec y$ is clearly a legal weighting for $H$, and $$\lambda'(H, \vec y)- \lambda'(H, \vec x)=x_j({\partial \lambda' (H, {\vec x})\over \partial x_i}- {\partial \lambda' (H, {\vec x})\over \partial x_j})-x_j^2{\partial^2 \lambda' (H, {\vec x})\over {\partial x_i \partial x_j}}=0.$$
So $\vec y$ is an optimal vector and $\vert\{i : y_i > 0 \}\vert= k-1$, contradicting the minimality of $k$. Hence Lemma \ref{Lemma2} holds. \qed

\bigskip
\noindent{\em  Proof of Theorem \ref{th1}.} Clearly, $\lambda'(H)\ge \lambda'(K^{\{1,2\}}_t)=2 - {1 \over t}$.

Now we proceed to show that $\lambda'(H)\le \lambda'(K^{ \{ 1,2 \} }_t)=2 - {1 \over t}$. Let ${\vec x}=(x_1, x_2, \ldots, x_n)$ be an optimal weighting of $H$ satisfying (*) with $k$ positive weights. Without loss of generality, we may assume that $x_1 \ge x_2 \ge \ldots \ge x_k >x_{k+1}=x_{k+2}=\ldots x_n=0$. By Lemma \ref{Lemma2}, $\forall 1\le i<j\le k,ij\in H^{2}$.

\begin{claim}\label{claim1}
$\forall 1\le i <j\le k$, if $i\in H$ but $j\notin H$, then $x_i-x_j=0.5$.
\end{claim}
 
\noindent{\em Proof of Claim \ref{claim1}.} By Lemma \ref{Lemma1}, ${\partial\lambda' (H, {\vec x})\over \partial x_i}={\partial \lambda' (H, {\vec x})\over \partial x_j}$. By Lemma \ref{Lemma2}, $\forall 1\le i<j\le k,ij\in H^{2}$, therefore $1+2(1-x_i)=2(1-x_j)$, i.e. $x_i-x_j=0.5$.  \qed

\begin{claim}\label{claim2}
Either $i\in H^{1}$ for all $1\le i \le k$ or  $i\notin  H^{1}$ for  all $1\le i \le k$.
\end{claim}

\noindent{\em Proof of Claim \ref{claim2}.} Assume that there are $l$ 1-sets of $\{1, 2, 3, \cdots, k\}$ in $H^1$. If $l=k$, then $i\in E^1$ for all $1 \le i\le k$, then $K^{\{1,2\}}_k$ is a subgraph of $H$. Since $t$ is the order of the maximum complete $\{1,2\}$-graph of $H$, then $k\le t$. We have
$$\lambda' (H, {\vec x})= \lambda' (K^{\{1,2\}}_k) = 2-{1\over k} \le 2-{1\over t}.$$

Therefore we can assume that  $l\le k-1$. Without loss of generality, assume that   $i \in H^1$ for $1\le i \le l $ and $i \notin H^1$ for $l+1\le j\le k$,  By Claim \ref{claim1}, $x_i=x_j+0.5$, $\forall 1\le i\le l$ and $l+1\le j\le k$. Then $l\le 1$. Otherwise, $x_1=x_k+0.5$ and $x_2=x_k+0.5$, contradicts to $\sum_{i=1}^{k} x_i =1$ and $ x_i
> 0 {\rm \ for \ } 1\le i\le k$. If $l=1$, then
$x_1=0.5+{0.5\over k}, x_2=x_3=\ldots =x_k={0.5\over k}$ and
\begin{eqnarray*}
\lambda' (H, {\vec x})&=& x_1+2\sum_{1\le i<j\le k}x_ix_j\\&=&0.5+{0.5\over k} + 2(0.5+{0.5\over k})(0.5-{0.5\over k})+2 {k-1\choose 2} ({0.5\over k})^2\\& =& 1.25+ {0.25\over k}- {0.5\over k^2}<1.5\\&\le& 2-{1\over t}.
\end{eqnarray*}
So Claim \ref{claim2} holds. \qed

Let's continue the proof of Theorem \ref{th1}.

If $i\in E^1$ for all $1 \le i\le k$, then $K^{\{1,2\}}_k$ is a subgraph of $H$. Since $t$ is the order of the maximum complete $\{1,2\}$-graph of $H$, then $k\le t$. We have
$$\lambda' (H, {\vec x})= \lambda' (K^{\{1,2\}}_k) = 2-{1\over k} \le 2-{1\over t}.$$

If $i\notin E^1$ for all $1 \le i\le k$, then $$\lambda' (H, {\vec x})= \lambda' (K^{(2)}_k) = 1-{1\over k} \le 2-{1\over t}.$$\qed

\subsection{Proof of Theorem \ref{th2}}

Let $F$ and $G$ be hypergraphs. We say that a function $f:V(F)\rightarrow V(G)$ is a $homomorphism$ from $F$ to $G$ if it preserves edges, i.e. $f(i_1)f(i_2)\cdots f(i_k)\in E(G)$ for all $i_1i_2\cdots i_k\in E(F)$. We say that $G$ is $F-hom-free$ if there is no homomorphism from $F$ to $G$. 
\begin{remark}\label{remarkhomfree1}
If $G$ is $F$-hom-free, then $G$ is $F$-free. 
\end{remark}

\noindent{\em Proof of Remark \ref{remarkhomfree1}.} If $G$ is  not $F$-free, then $G$ contains a copy of $F$  as a subgraph. Let $f:V(F)\rightarrow V(G)$ be the function defined by $f(v)=v$ for every $v\in V(F)$. Then $f$ is a $homomorphism$ from $F$ to $G$. So $G$ is not $F$-hom-free. \qed

\begin{remark}\label{remarkhomfree2}
$G$ is $F$-hom-free if and only if the blowup $G(s, s, \cdots, s)$ is $F$-free for every $s$. 
\end{remark}

\noindent{\em Proof of Remark \ref{remarkhomfree2}.} If $G$ is  not $F$-hom-free, then there exists a function $f:V(F)\rightarrow V(G)$ which is a homomorphism from $F$ to $G$. Let $s=\max\{\vert f^{-1}(v)\vert, v\in V(G)\}$. Then  $G(s, s, \cdots, s)$ contains $F$ as a subgraph. 

  Assume that  $G(s, s, \cdots, s)$ contains $F$ as a subgraph for some $s$. Then for each $v\in V(F)$, $v$ is contained in a set of some vertices of $G(s, s, \cdots, s)$ blowed up by a vertex $w \in V(G)$. Let  $f(v)=w$. Then $f$ is a homomorphism from $F$ to $G$. \qed

\begin{remark}\label{completehomfree} If $H$ is a $\{1,2\}$-graph and $t=\chi(H^2)$, then a complete $\{1,2\}$-graph  $K_l^{\{1,2\}}$ is $H$-hom-free if and only if $l\leq t-1$.
\end{remark} 

\noindent{\em Proof of Remark \ref{completehomfree}.} Apply Remark \ref{remarkhomfree2}. \qed

\medskip

We can make an analogous definition to the Tur\'{a}n density:
 \begin{eqnarray*}
 \pi_{hom}(F)=\lim_{n\rightarrow\infty}\max \{h_n(G): |v(G)|=n, G\subseteq K_{n}^{R(F)}, {\ \rm and \ } G {\rm \ is\ } F{\rm \ -hom-free.}\}.
 \end{eqnarray*}

Then we have two useful lemmas.

\begin{lemma}\label{Lemma3}
$\pi_{hom}(F)=\pi(F).$
\end{lemma}

\noindent{\em Proof of Lemma \ref{Lemma3}.} Let $R(F)=R$ and $G\subseteq K_{n}^{R}$. If $G$ is $F$-hom-free, then by Remark \ref{remarkhomfree1}, $G$ is $F$-free. So $\pi(F) \ge \pi_{hom}(F)$. On the other hand,
$\forall \varepsilon >0,\forall n_0,\exists n>n_0,\exists G \subseteq K_n^R$  and $G$ is not $F$-hom-free such that $h_n(G)\leq \pi_{hom}(F)+\varepsilon $. Since there is a homomorphism from $F$ to $G$ with $n$ vertices, then by Remark \ref{remarkhomfree2}, there exists $s$ such that $G(s,s,\ldots,s)$ contains $F$. So,
\begin{eqnarray*}
\pi(F)&\leq& \lim_{s\rightarrow\infty}h_{ns}(G(s,s,\ldots,s))\\&= &h_n(G)\\&\leq&\pi_{hom}(F)+\varepsilon.
\end{eqnarray*}
Hence, $\pi(F)\leq \pi_{hom}(F).$\qed

\begin{lemma}\label{Lemma4}
$\pi(F)$ is the supremum of $\lambda'(G)$ over all $F$-hom-free $G$ with $R(G)\subseteq R(F)$.
\end{lemma}

\noindent{\em Proof of Lemma \ref{Lemma4}.} Suppose that $F$ is a hypergraph and $G$ is an $F$-hom-free hypergraph with $n$ vertices and $R(G)\subseteq R(F)$. Let ${\vec s}=(s_1,s_2,\ldots,s_n)$ be an optimal vector of $\lambda'(G)$. Take any $m$, note that $G(s_1m,s_2m,
\ldots,s_nm)$ is an $F$-free hypergraph on $p=({\sum_{i=1}^n{s_i}})m$ vertices with $R(G)\subseteq R(F)$ and 
$$h_{p}(G(ms_1,ms_2,\ldots,ms_m))=\sum_{j\in R(G)}\sum_{i_1i_2\ldots i_j\in E(G)}{{(s_{i_1}s_{i_2}\ldots s_{i_j})m^j}\over{{sm \choose j}}}\rightarrow \lambda'(G,{\vec s}) {\rm \ as\ } m\rightarrow \infty.$$
So $\pi(F)\geq \lambda'(G,\vec{s})=\lambda'(G)$.

On the other hand, $\forall \varepsilon>0$, $\forall n_0$, $\exists n>n_0$, $\exists$ an $F$-hom-free $H$ with $n$ vertices and $R(H)\subseteq R(F)$ such that $\pi(F)\leq h_n(H) +\varepsilon$. Note that
\begin{eqnarray*}
\lambda'(H)&\ge& \lambda'(H,({1\over n},{1\over n},\ldots,{1\over n}))\\&=&\sum_{j\in R(H)}j! {\sum_{i_1i_2\ldots i_j\in E(H)}({1\over n})^j}\\&\ge& h_n(H)-\varepsilon {\rm \ when\ }n {\rm \ is \ \ large\ \ enough}.
\end{eqnarray*}
So, $\pi(F)\leq \lambda'(H) +2\varepsilon$.

Therefore, $\pi(F)$ is the supremum of $\lambda'(G)$ over all $F$-hom-free hypergraphs $G$ with $R(G)\subseteq R(F)$ .\qed

To continue the proof of Theorem \ref{th2}, we define a dense hypergraph. 
 
\begin{defi} \label{2.1}
A hypergraph $G$ is $dense$ if every proper subgraph $G'$ satisfies $\lambda'(G')<\lambda'(G)$.
\end{defi}

\begin{remark}\label{re1}By Theorem \ref{MStheo} and Remark \ref{relation}, a graph $G$ is dense if and only if $G$ is $K_t^{\{2\}}$. By Theorem \ref{th1}, a $\{1,2\}$-graph $G$ is dense if and only if $G$ is $K_t^{\{1,2\}} \,\  (where \,\ t\geq 2)$.
\end{remark}

 \noindent{\em Proof of Theorem \ref{th2}.} Assume that $H$ is a $\{1,2\}$-graph and $H^2$ is not bipartite. By lemma \ref{Lemma4}, $\pi(H)$ is the supremum of the Lagrangians of  all $H$-hom-free $\{1,2\}$-graphs, all $H$-hom-free graphs and all $H$-hom-free $\{1\}$-graphs . So $\pi(H)$ is the supremum of the Lagrangians of  all dense $H$-hom-free  $\{1,2\}$-graphs, all dense $H$-hom-free   graphs and all $H$-hom-free $\{1\}$-graphs. Let $t=\chi(H^2)\ge 3$. By Remark \ref{re1} and \ref{completehomfree}, a dense $H$-hom-free  $\{1,2\}$-graph must be $K_l^{\{1,2\}},2\le l\leq t-1$ and a dense $H$-hom-free  graph must be $K_s$. Also, note that the Lagrangian of all $\{1\}$-graphs is 1. So,
 $$\pi(H)=\max\{\lambda'(K_{t-1}^{\{1,2\}}), \lambda' (K_s), 1 \}=\max\{2-{1\over {t-1}}, 1-{1 \over s}, 1\}=2-{1\over{t-1}}.$$\qed

\bigskip
{\bf Acknowledgments.} This research is partially  supported by National Natural Science Foundation of China (No. 11271116).

\end{document}